\documentclass[12pt]{article}
\usepackage{diagrams}
\newtheorem{proposition}{Proposition}
\newtheorem{theorem}[proposition]{Theorem}
\newtheorem{lemma}[proposition]{Lemma}
\newtheorem{corollary}[proposition]{Corollary}

\newtheorem{remark}[proposition]{Remark}

\newcommand{\dem}{\par\noindent{\it Proof. }}

\newcommand{\spa}{\vspace{2.0mm}}
\begin{document}

\title{Regularity and finite injective dimension in characteristic $p>0$}
\author{Tiberiu Dumitrescu and Cristodor Ionescu\footnote{This paper was written while the authors were visiting the \textit{School of Mathematical Sciences}, Lahore. They thank this institution for support and hospitality. }\\[2mm]
 Faculty of Mathematics, University of Bucharest, \\14 
Academiei St., Bucharest, RO-010014, Romania\\
email tiberiu@al.math.unibuc.ro\\[2mm]
Institute of Mathematics Simion Stoilow of the Romanian Academy,\\
P.O. Box 1-764, Bucharest, RO-014700, Romania
\\
email Cristodor.Ionescu@imar.ro
}
\date{}

\maketitle
\begin{abstract}
Recently,  the regular   local rings of prime characteristic were cha\-rac\-te\-rized in
terms of the finiteness of injective dimension of the Fro\-be\-ni\-us map. We obtain relative versions of
this result.
\end{abstract}

Let $A$ be a noetherian local ring of prime characteristic $p>0$ and   let $F_A$ denote the Frobenius morphism of $A$. A celebrated theorem of Kunz \cite{Ku} asserts that $A$ is regular if and only if $F_A$ is flat. Later, Rodicio \cite{Ro} has shown that  $A$ is regular if and only if $F_A$ has finite flat dimension. Recently Avramov, Iyengar and Miller \cite{AvIyMi} have shown that $A$ is regular if and only if $F_A$ has finite injective dimension.
\par Andr\'e and Radu \cite{An1},\cite{Ra} proved  a relative version of Kunz's result. In order to state it, we consider the following setup. For a ring $C$ of prime characteristic $p>0,$  denote by $C^{(p)}$ the $C$-algebra structure on $C$ given by $F_C.$ Now let  $u:A\rightarrow B$ be a morphism of noetherian rings of prime characteristic $p>0.$ 
We get the  commutative diagram

\begin{diagram}[size=2.5em]
 A       & \rTo<{u}   &B &            &        &     &       \\
 \dTo<{F_A}     &        &  \dTo>{F_A\otimes_A B}              &            &        &     &       \\
A^{(p)}   &\rTo<{u\otimes_A A^{(p)}}    & A^{(p)}\otimes_A B &\rTo^{w_{B/A}}  & B^{(p)}& & 
\end{diagram}
where $w_{B/A}$ is given by $w_{B/A}(a\otimes b)=ab^p,$ for every $a\in A,\ b\in B.$
Note that $w_{B/A}$ induces a bijection on the spectra. In particular, if $w_{B/A}$ is flat and $B^{(p)}$ is noetherian, then $A^{(p)}\otimes_A B$ is also noetherian.
The following result was proved by Andr\'e,   Radu   and the first author of the present paper. 
See \cite{An1} and \cite{Ra} for the equivalence of $(a)$ and $(b)$, and \cite[Th. 2.13]{Du2} for $(c)$.

\begin{theorem}\label{1.0.1}
Let  $u:A\rightarrow B$ be a  morphism of noetherian  rings. The following assertions are
equivalent: 

$(a)$ u is regular;

$(b)$  $w_{B/A}$ is flat;

$(c)$ u is flat and   $w_{B/A}$ has finite flat dimension.
\end{theorem}

The aim of this note is to investigate  relative versions of the result of Avramov, Iyengar and Miller described above, in the spirit of Theorem \ref{1.0.1}.

\spa
All rings considered in this note have prime characteristic $p$. Notations and terminology are as in
\cite{Mat}. Consider the above setup. If the ring $A^{(p)}\otimes_A B$ is  noetherian,
from the previous results we obtain:

\begin{theorem}\label{1}
Let  $u:A\rightarrow B$ be a flat local morphism of noetherian local rings. Suppose that $A^{(p)}\otimes_A B$ is  noetherian.  Then   $w_{B/A}$ has finite injective dimension if and only if u is regular and $A$ is Gorenstein.   
\end{theorem}

\dem If $w_{B/A}$ has finite injective dimension,   \cite[Th. 13.2]{AvIyMi} implies that $A^{(p)}\otimes_A B$ is Gorenstein and $w_{B/A}$ has finite flat dimension. Applying Theorem \ref{1.0.1}, we obtain that $u$ is regular. Since $u$ is faithfully flat, $A^{(p)}\otimes_A B$ is also faithfully flat over  $A^{(p)}.$ Hence $A$ is Gorenstein.

Conversely, since $u$ is regular and $A$ is Gorenstein, $B$ is also Gorenstein. By Theorem \ref{1.0.1} it follows that $w_{B/A}$ is faithfully flat, hence $A^{(p)}\otimes_A B$ is Gorenstein. \cite[Th. 13.2]{AvIyMi} shows that $w_{B/A}$ has finite injective dimension.$\bullet$

\begin{corollary}\label{1.1}
Let A be a noetherian local ring with geometrically regular formal fibers. Then   $A$ is Gorenstein if and only if $\widehat{A}$ has finite injective dimension over $A[\widehat{A}^p]$.   
\end{corollary}

\dem From the assumption and Theorem \ref{1.0.1} it folows that $w_{\widehat{A}/A}$ is flat, hence injective. 
Then the image of $w_{\widehat{A}/A}$ is $A[\widehat{A}^p].$ Now apply Theorem \ref{1}.$\bullet$

\begin{remark}\label{2.0}
{\rm Consider the diagram before Theorem \ref{1.0.1}. Let $A'$ be an $A$-algebra and set $B'=B\otimes_A A'.$ Then it is easy to check that $w_{B'/A'}=w_{B/A}\otimes_{A^{(p)}} A'^{(p)}$.}
\end{remark}

\spa  

According to \cite[Def. 2.8]{Du2} a local ring is called \textit{quasi-basic} if the $A$-algebra $A^{(p)}$ is a directed union of a family $(A_i)_{i\in I}$ of  finite subalgebras, such that $A^{(p)}$ is flat over $A_i,$ for any $i\in I.$ By \cite{Du2}, any  ring with Artin's approximation property (e.g. a complete local ring) is  quasi-basic. Also, by \cite[Prop. 2.9]{Du2}, any algebra essentially of finite type over a quasi-basic ring is quasi-basic. In particular, algebras which are essentially of finite type over a field are quasi-basic. Denote by ${\rm fd}_A(M)$ (resp. ${\rm id}_A(M)$), the flat (resp. injective) dimension of an $A$-module $M.$ 

\begin{lemma}\label{100}
Let $f:R\rightarrow S$ be a local morphism of local noetherian rings such that ${\rm id}_R(S)<\infty$. Then 
${\rm fd}_R( S)\leq{\rm id}_R(S)+{\rm edim}(S).$
\end{lemma}

\dem Consider a minimal Cohen factorization of the canonical morphism $R\rightarrow\widehat{S}$ (see \cite{AvIyMi})

\begin{diagram}[size=2.0em]
 R      & \rTo   &S &            &        &     &       \\
 \dTo     &        &  \dTo              &            &        &     &       \\
R'   &\rTo    & \widehat{S}. & & & & 
\end{diagram}
From \cite[Cor. 2.5]{AvIyMi} we get 
$${\rm id}_{R'}(\widehat{S})\leq{\rm id}_R(S)+{\rm edim}(S)<\infty.$$
By \cite[Th. 13.2]{AvIyMi}, $R'$ is a Gorenstein ring. From \cite[Th. 2.2]{LeVa} we obtain
${\rm fd}_{R'}(\widehat{S})<\infty$ so, by Auslander-Buchsbaum formula, 
$${\rm fd}_{R'}(\widehat{S})\leq{\rm depth}(R')={\rm id}_{R'}(\widehat{S}).$$ So, from \cite[Cor. 2.5]{AvIyMi} we obtain 
$${\rm fd}_{R}(S)\leq{\rm fd}_{R'}(\widehat{S})\leq{\rm id}_{R'}(\widehat{S})\leq{\rm id}_{R}(S)+{\rm edim}(S).\bullet$$

\begin{theorem}\label{2}
Let  $u:A\rightarrow B$ be a local morphism of noetherian local rings. Assume  that A is   quasi-basic and  $w_{B/A}$ has finite injective dimension. Then $w_{B/A}$ has finite flat dimension.  Moreover if  $u$ is flat, then $u$ is  regular.
\end{theorem}

\dem Set  $w:=w_{B/A}$ and $C:= A^{(p)}\otimes_A B.$ Since A is a  quasi-basic ring, the $A$-algebra $A^{(p)}$ is a directed union of a family $(A_i)_{i\in I}$ of  finite subalgebras such that $A^{(p)}$ is flat over $A_i,$ for any $i\in I.$ Then $B_i:=A_i\otimes_A B$ is a noetherian local ring. 
Consider the following commutative diagram with canonical morphisms:

\begin{diagram}[size=1.6em]
  A       & \rTo^u  & B                  &            &        &     &       \\
 \dTo<{\tau_i}     &        &\dTo>{v_i}                &            &        &     &       \\
A_i       & \rTo   &B_i &            &        &     &       \\
 \dTo<{\sigma_i}     &        &  \dTo>{s_i}              &            &        &     &       \\
A^{(p)}   &\rTo    & C&\rTo^w  & B^{(p)}.& & 
\end{diagram} 
Since $\sigma_i$ is flat, it follows that $s_i=\sigma_i\otimes_A B$  is  faithfully flat.  Hence  
$${\rm Ext}_C^n(E\otimes_{B_i}C,B^{(p)})\cong{\rm Ext}_{B_i}^n(E,B^{(p)})$$
for any $B_i$-module $E$ and for any $n\in{\bf N},$ according to \cite[Th. 11.65]{Rot}. 
It follows that ${\rm id}(ws_i)\leq {\rm id}(w).$
From Lemma \ref{100} we have 
$${\rm fd}(ws_i)\leq{\rm id}(ws_i)+{\rm edim}(B):=k.$$
 Let 
$j>k.$  Then ${\rm Tor}_j^{B_i}(B^{(p)},E)=0,$ for any $B_i$-module $E.$ A direct limit argument (see \cite[Ch. VI, Exercise 17]{CaEi}) shows that 
${\rm Tor}_j^C(B^{(p)},F)=0,$ for any $C$-module $F.$ Thus ${\rm fd}(w)\leq k.$
 In addition, if $u$ is flat, then $u$ is regular by  Theorem \ref{1.0.1}.$\bullet$

\begin{corollary}\label{2.1}
Let  $u:A\rightarrow B$ be a  local morphism of noetherian local rings, such that A is a  quasi-basic regular ring.  If  $w_{B/A}$ has finite injective dimension, then $u$ is a regular morphism.    
\end{corollary}

\dem Apply Theorem \ref{2} and \cite[Th. 2.10]{Du2}.$\bullet$

\spa

Let $u:A\rightarrow B$ be a  morphism of noetherian local rings. It is well-known that  when $A$ is a field, $u$ is a regular morphism if and only if $B$ is a formally smooth $A$-algebra.  From Theorems \ref{1} and \ref{2} we get the following consequence:

\begin{corollary}\label{3}
Let  k be a field and B  a  noetherian local k-algebra. Then B is a formally smooth k-algebra if and only  if $w_{B/k}$ has finite injective dimension.$\bullet$    
\end{corollary}
 
 When $k={\bf Z}_p,$ it follows that $w_{B/k}=F_B,$ hence we retrieve \cite[Th. 13.3]{AvIyMi}: A noetherian local ring B is  regular  if and only if $F_B$ has finite injective dimension.

\end{document}